\documentclass[letterpaper, 10 pt, conference]{IEEEconf}
\IEEEoverridecommandlockouts

\usepackage{cite}
\usepackage{amsmath,amssymb,amsfonts}

\usepackage{accents}
\usepackage{algorithmic}
\usepackage{graphicx}

\usepackage{subfigure}
\usepackage{epsfig}
\usepackage{cite}

\usepackage[hidelinks]{hyperref}

\usepackage[dvipsnames]{xcolor}
\usepackage{tikz}
\usepgflibrary{fpu}
\usetikzlibrary{positioning,calc,arrows,shapes,shadows,matrix,backgrounds}
\usetikzlibrary{external}

\usetikzlibrary{positioning,calc,arrows,shapes,shadows}
\tikzset{
	auto,
	node/.style={
		circle,
		draw,
		fill=gray!25,
		minimum size = 1.2*\dh,
		inner sep=\dn},
	sum/.style={circle,draw,draw=black,inner sep=0mm,minimum size=2mm},
	jun/.style={circle,draw,draw=black,inner sep=0mm,minimum size=0mm},
	>={latex},
}
\def\dn{1ex*0.75}

\def\dh{6*\dn}

\usepackage{pgfplots}
\pgfplotsset{compat=1.18}

\usepackage{mathtools}

\newcommand{\set}[1]{\mathbb{#1}}
\newcommand{\N}{\set{N}}
\newcommand{\R}{\set{R}}

\newcommand{\SI}{\set{I}}
\newcommand{\SU}{\set{U}}
\newcommand{\SX}{\set{X}}

\newcommand{\CA}{\mathcal A}
\newcommand{\CB}{\mathcal B}

\newcommand{\CD}{\mathcal D}

\newcommand{\CK}{\mathcal K}
\newcommand{\CL}{\mathcal L}

\newcommand{\CS}{\mathcal S}

\newcommand{\CX}{\mathcal X}

\newcommand{\eps}{\varepsilon}
\newcommand{\lb}{\lambda}
\newcommand{\al}{\alpha}

\newcommand{\de}{\delta}
\newcommand{\sg}{\sigma}

\newcommand{\abs}[1]{\left| #1 \right|}
\newcommand{\floor}[1]{\lfloor #1 \rfloor}
\newcommand{\norm}[1]{\mbox{$\left\| #1 \right\|$}}

\newtheorem{theo}{Theorem}
\newtheorem{prop}[theo]{Proposition}
\newtheorem{hypo}[theo]{Hypothesis}
\newtheorem{lemm}[theo]{Lemma}
\newtheorem{conj}[theo]{Conjecture}
\newtheorem{coro}[theo]{Corollary}
\newtheorem{defi}[theo]{Definition}
\newtheorem{rema}[theo]{Remark}
\newtheorem{assu}[theo]{Assumption}
\newtheorem{exam}{Example}
\newcommand{\theorem}[2][]{\begin{theo}[#1]#2\end{theo}}
\newcommand{\proposition}[2][]{\begin{prop}[#1]#2\end{prop}}

\newcommand{\lemma}[2][]{\begin{lemm}[#1]#2\end{lemm}}

\newcommand{\definition}[2][]{\begin{defi}[#1]#2\end{defi}}

\newcommand{\ass}[2][]{\begin{assu}[#1]#2\end{assu}}
\newcommand{\example}[1]{\begin{exam}#1\end{exam}}

\renewcommand{\t}{\tilde}
\newcommand{\h}{\hat}
\newcommand{\ubar}[1]{\underline{#1}}
\newcommand{\obar}[1]{\overline{#1}}

\newcommand\copyrighttext{%
	\footnotesize 
	\textcopyright 2026 IEEE. Personal use of this material is permitted.
	Permission from IEEE must be obtained for all other uses, in any current or future media, including reprinting/republishing this material for advertising or promotional purposes, creating new collective works, for resale or redistribution to servers or lists, or reuse of any copyrighted component of this work in other works.
}
\newcommand\copyrightnotice{%
	\begin{tikzpicture}[remember picture,overlay]
		\node[anchor=south,yshift=10pt] at (current page.south) {\fbox{\parbox{\dimexpr\textwidth-\fboxsep-\fboxrule\relax}{\copyrighttext}}};
	\end{tikzpicture}%
}

\begin{document}

	\title{On stability and non-averaged performance of economic MPC with terminal conditions for optimal periodic operation}

	\author{
		Jonas Mair$^{a}$,
		Lukas Schwenkel$^{a}$,
		Matthias A. M\"uller$^{b}$,
		Frank Allg\"ower$^{a}$
		\thanks{$^{a}$ University of Stuttgart, Institute for Systems Theory and Automatic Control, 70550 Stuttgart, Germany; (e-mail: {\it jonas.mair@ist.uni-stuttgart.de, lukas.schwenkel@ist.uni-stuttgart.de, frank.allgower@ist.uni-stuttgart.de})}
		\thanks{$^{b}$  Leibniz University Hannover, Institute of Automatic Control, 30167 Hannover, Germany; (e-mail: {\it mueller@irt.uni-hannover.de})}
	}

	\maketitle
	\copyrightnotice
	\thispagestyle{empty}

	\begin{abstract}
		Operation at steady state is often not optimal when optimizing over an economic cost objective. In many cases, periodic operation yields better performance. Therefore, we derive asymptotic stability guarantees of an economic model predictive control scheme with terminal conditions for systems with optimal periodic operation for a more general setup than existing methods can handle. Moreover, we establish a non-averaged closed-loop performance bound by defining the closed-loop cost via a Cesàro summation instead of ordinary summation. Such a non-averaged performance bound provides new insights for systems with periodic optimal operation.
	\end{abstract}
	
	\section{Introduction}
	
	Economic Model Predictive Control (EMPC) (see e.g. \cite{Ellis2014},\cite{Faulwasser2018}) is well understood for systems with optimal operation at steady state. For this case, asymptotic stability and practical asymptotic stability guarantees exist for EMPC schemes with and without terminal conditions \cite{Amrit2011}, \cite{Gruene2013}. Moreover, conditions for optimal and approximately optimal averaged performance were established with and without terminal conditions \cite{Angeli2012}, \cite{Gruene2014} and the same holds for non-averaged performance results \cite{Gruene2015}, \cite{Gruene2014}. However, even if an EMPC scheme achieves optimal averaged performance, the transient costs can still be arbitrarily bad on any finite time horizon. Thus, non-averaged performance results are of interest. The derivation of such results often relies on a strict dissipativity assumption, which is also tightly connected to the existence of an optimal steady state~\cite{Mueller2015}. When optimizing over economic costs, steady state operation is, however, not necessarily optimal. Therefore, \cite{Mueller2016,Schwenkel2024a,Schwenkel2024,Zanon2016,Dong2018} consider the more general case of optimal periodic operation. Also for this case, there exist EMPC schemes with stability and performance guarantees. With a multi-step EMPC scheme without terminal conditions, optimal averaged performance and convergence to the optimal periodic orbit can be guaranteed, however, it is not necessarily stabilized \cite{Mueller2016}. Instead, one can use discounted cost formulations to achieve practical asymptotic stability \cite{Schwenkel2024a}. A special case of these discounts are linear discounts which achieve the same guarantees \cite{Schwenkel2024} and are equivalent to a so-called Cesàro cost \cite{Mair25}. Moreover, with terminal conditions, one can asymptotically stabilize a fixed phase of a periodic orbit \cite{Zanon2016}, or obtain asymptotic stability and averaged performance guarantees for even more general operating behaviors \cite{Dong2018}. In \cite{Dong2018}, periodic orbits are also examined as a special case. However, besides a strict dissipativity assumption, this work requires the control law to be continuous in order to show asymptotic stability, which can be restrictive and is not necessarily needed for the steady state case. Furthermore, all the previous mentioned results rely on the optimal periodic orbit to be minimal (cf. Definition~\ref{def:periodic_orbits} below). While this is often an explicit assumption, in \cite{Dong2018} the continuity of the optimal control law implicitly excludes non-minimal periodic orbits. Furthermore, non-averaged performance guarantees are not available for systems with periodic optimal operation.
	
	In this work we derive asymptotic stability guarantees of an EMPC scheme with terminal conditions for systems which are optimally operated at a periodic orbit, without requiring minimality of the optimal periodic orbit. We also introduce a suitable definition of asymptotic stability for this case. Moreover, for minimal optimal periodic orbits, we establish infinite-horizon non-averaged performance bounds by defining the closed loop-cost via Cesàro summation instead of ordinary summation, which is argued to be a suitable measure and a natural extension for periodic operation \cite{Mair25}. At last, we show that an infinite-horizon Cesàro cost is a suitable choice for the terminal cost and therefore provide an understanding how such terminal costs can be designed.

	\emph{Outline:} In Section~\ref{sec:MPC-scheme}, the EMPC scheme is introduced, before deriving asymptotic stability guarantees in Section~\ref{sec:stability} and infinite-horizon non-averaged performance bounds in Section~\ref{sec:performance}. The connection of Cesàro costs and suitable terminal costs is discussed in Section~\ref{sec_cesaro_terminal}. The paper is concluded in Section~\ref{sec:conclusion}.
	
	\emph{Notation:} For $a\leq b$, the set of integers in the interval $[a,b]$ is denoted $\SI_{[a,b]}$. We denote the modulo operation by $[k]_p$, i.e., the remainder when dividing $k$ by $p$. For $x\in\R$, the floor operator $\floor{x}$ crops all decimal places. For continuous functions $\al,\de:[0,\infty)\to[0,\infty)$, we say that $\al\in\CK_\infty$ if and only if $\al(0)=0$, $\lim_{t\rightarrow\infty}\al(t)=\infty$ and $\al$ is strictly increasing and we say $\de\in\CL$ if and only if $\de$ is monotonically decreasing and $\lim_{t\rightarrow\infty}\de(t)=0$. For a function $h:\CA\to\CB$, the image of $h$ on $\CD\subseteq\CA$ is defined as $h(\CD)\coloneq\{h(d)\,|\,d\in\CD\}$. 
	The cardinality of a set $\SX$ is measured by $\#\SX$.
	

	\section{EMPC-scheme with terminal cost} \label{sec:MPC-scheme} 
	Consider the nonlinear discrete-time system
	\begin{gather}\label{eq:sys}
		x(k+1)=f(x(k),u(k)),
	\end{gather}
	subject to state and input constraints $x(k)\in\SX\subset\R^n$ and $u(k)\in\SU\subset\R^m$. A state trajectory of \eqref{eq:sys} resulting from an initial condition $x_0\in\SX$ and a specific input sequence $u\in\SU^N$ of length $N\in\N$ is denoted $x_u(k,x_0)$ and defined by $x_u(0,x_0)=x_0$ and $x_u(k+1,x_0)=f(x_u(k,x_0),u(k))$ for $k\in\SI_{[0,N-1]}$. The set of all feasible control sequences of length $N$ and starting at $x\in\SX$ is denoted $\SU^N(x)\coloneq\{u\in\SU^N~|~\forall k\in\SI_{[0,N]}:x_u(k,x)\in\SX\}$. For a given stage cost $\ell:\SX\times\SU\rightarrow\R$ and a terminal cost  $V_\mathrm{f}(x):\SX_\mathrm{f}\to\R$, which is defined on a terminal constraint set $\SX_\mathrm{f}\subseteq\SX$, we define the finite horizon cost functional
	\begin{align} \label{eq:finite_horion_cost_functional}
		J_N(x,u)=\sum_{k=0}^{N-1} \ell(x_u(k,x),u(k))
		+ V_\mathrm{f}(x_u(N,x)),
	\end{align}
	The corresponding optimal value function is defined by
	\begin{align}\label{eq:finite_horizon_value_function}
		V_N(x)\coloneq\inf_{\substack{u\in\SU^N(x) \\ x_u(N,x)\in\SX_\mathrm{f}}} J_N(x,u)
	\end{align}
	and the set of feasible states is denoted $\SX_N\coloneq\{x\in\SX\,|\,\exists u\in\SU^N(x):x_u(N,x)\in\SX_\mathrm{f}\}$. 
	
	Suppose for a fixed horizon $N\in\N$ and any $x\in\SX_N$ there exists an optimal input sequence $u^\star_{N,x}\in\SU^N(x)$ satisfying\footnote{Existence of such an optimal input sequence is, e.g., guaranteed under Assumption~\ref{ass:cont_comp} below.} $J(x,u^\star_{N,x})=V_N(x)$. Then, using the finite horizon cost functional~\eqref{eq:finite_horion_cost_functional} and the optimal value function~\eqref{eq:finite_horizon_value_function}, an EMPC scheme can be defined as follows. Measure the current state $x=x(n)$ of the system. Then, minimize $J_N(x,u)$ over $u\in\SU^N(x)$ to obtain the optimal value function and denote the minimizing input sequence $u^\star_{N,x}$. Finally, the first element of $u^\star_{N,x}$ is applied until the next time instant. This leads to the EMPC feedback law
	\begin{align} \label{eq:feedback_law}
		\mu_N(x)\coloneq u^\star_{N,x}(0),
	\end{align}
	i.e., the closed-loop system can be written as $x(n+1)=f(x_{\mu_N}(n,x),\mu_N(x_{\mu_N}(n,x)))$. Note that the optimal feedback law $\mu_N(x)$ is in general neither unique nor continuous, since the minimizer $u^\star_{N,x}$ is not necessarily unique.

	\section{Stability guarantees} \label{sec:stability}
	
	In this section, we derive an asymptotic stability result for the EMPC scheme in Section~\ref{sec:MPC-scheme} for systems which are optimally operated at a (possibly non-minimal) periodic orbit, by suitably extending the arguments and assumptions in \cite{Amrit2011}. For minimal optimal periodic orbits, this issue has been addressed in \cite{Zanon2016,Dong2018}. In \cite{Zanon2016}, a periodic strict dissipativity assumption with phase dependent storage function is employed. Furthermore, a phase-dependent terminal cost is used and therefore also asymptotic stability is only guaranteed w.r.t. a particular phase of the optimal periodic orbit. Hence, such an EMPC scheme may leave the optimal orbit temporarily to return in a different phase due to the artificial phase dependency in the terminal conditions. In contrast, \cite{Dong2018} introduces these quantities independent of the orbit's phase and therefore also the derived stability guarantees are independent of the orbit's phase, which is why we will mostly relate our following work to \cite{Dong2018} instead of \cite{Zanon2016}.

	\subsection{Assumptions} \label{sec:assumptions}
	
	First, we need continuity of the system dynamics, the stage costs and the terminal cost as well as compactness of the state-input constraint set.
	\ass[Continuity and compactness]{ \label{ass:cont_comp}
		The constraint set $\SX\times\SU$ is compact and  the system dynamics $f$ and the stage cost $\ell$ are continuous on $\SX\times\SU$. Furthermore, the terminal cost function $V_\mathrm{f}$ is continuous on $\SX_\mathrm{f}$.
	}

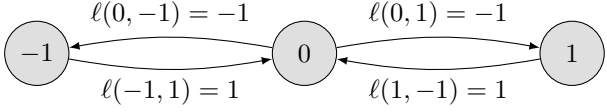
\begin{figure}\center
		\begin{tikzpicture}[xscale=1,yscale=1]
			\node[node] (n0) at (0,0)  {$0$};
			\node[node] (n1) at (-5*\dh,0)  {$-1$};
			\node[node] (n2) at (5*\dh,0)  {$1$};
			\draw[->] (n2) to [out= -170, in= -10] node[midway,below] {$\ell(1,-1)=1$} (n0);
			\draw[->] (n0) to [out= 10, in= 170] node[midway,above] {$\ell(0,1)=-1$} (n2);
			\draw[->] (n1) to [out=-10, in=-170] node[midway,below] {$\ell(-1,1)=1$} (n0);
			\draw[->] (n0) to [out= 170, in= 10] node[midway,above] {$\ell(0,-1)=-1$} (n1);
		\end{tikzpicture}
		\caption{
			Illustration of the non-minimal optimal periodic orbit $\varPi^\star$.
		}
		\label{fig:non_minimal_orbit}
	\end{figure}

	Since we are interested in systems with optimal operation at a periodic orbit, let us now formally define (optimal and minimal) periodic orbits analogously to \cite{Schwenkel2024a}.
	\definition[Periodic Orbit]{ \label{def:periodic_orbits}
		Let $\varPi_\SX:\SI_{[0,p-1]}\rightarrow\SX$ and $\varPi_\SU:\SI_{[0,p-1]}\rightarrow\SU$ for $p\in\N$. Then $\varPi=(\varPi_\SX,\varPi_\SU)$ is called a {\it$p$-periodic orbit} of system \eqref{eq:sys}, if 
		\begin{align}\label{eq:periodic_orbit}
			\varPi_\SX([i+1]_p)=f(\varPi_\SX(i),\varPi_\SU(i))=f(\varPi(i))
		\end{align}
		for all $i\in\SI_{[0,p-1]}$. A $p$-periodic orbit is called {\it minimal} if $\varPi_\SX$ is injective. The distance of a point $(x,u)\in{\SX\times\SU}$ to the $p$-periodic orbit $\varPi$ is defined as $\norm{(x,u)}_\varPi\coloneq\min_{i\in\SI_{[0,p-1]}} \norm{(x, u)-\varPi(i)}$. Furthermore, the closest point of $\varPi_\SX$ to a given point $x\in\SX$ is characterized by $i_x\in\SI_{[0,p-1]}$ satisfying $\norm{x-\varPi(i_x)}=\norm{x}_{\varPi_\SX}$. The set of all feasible $p$-periodic orbit is denoted $S_\varPi^p$. A feasible $p^\star$-periodic orbit $\varPi^\star$ is called an optimal periodic orbit, if
		\begin{gather*} \label{eq:opt_orbit}
			\ell_{\varPi^\star}\coloneq\frac{1}{p^\star}\sum_{i=0}^{p^\star-1}\ell(\varPi^\star(i))=\inf_{p\in\N,\varPi\in S_\varPi^p} \frac{1}{p}\sum_{i=0}^{p-1}\ell(\varPi(i)).
		\end{gather*}
	}
	\vspace*{\belowdisplayskip}
	To derive stability results in EMPC, often a strict dissipativity assumption is used. We use an extension of strict dissipativity for steady states to periodic orbits \cite[Ass. 1]{Kohler2018}.
	\ass[Strict dissipativity]{\label{ass:strict_dissipativity}
		There exist a continuous storage function $\lb:\SX\rightarrow\R$ and a function $\ubar{\al}_{\t\ell}\in\CK_\infty$, such that the rotated stage cost $\t\ell(x,u)\coloneq\ell(x,u)-\ell_{\varPi^\star}+\lb(x)-\lb(f(x,u))$ satisfies for all $x\in\SX$ and $u\in\SU^1(x)$
		\begin{equation}\label{eq:rotated_cost_pos}
			\t\ell(x,u)\geq\ubar{\al}_{\t\ell}(\norm{(x,u)}_{\varPi^\star}).
		\end{equation}
	}
	\vspace*{\belowdisplayskip}
	Strict dissipativity is a sufficient condition for optimal operation at a periodic orbit, i.e., that the optimal asymptotic average performance is $\ell_{\varPi^\star}$ and under an additional local controllability assumption it is also necessary \cite{Mueller2015}. Note that the argument of $\al_{\t\ell}$ depends on $u$ and therefore our strict dissipativity assumption is stronger than the one used in \cite{Dong2018}.
	To guarantee asymptotic stability, we furthermore need a similiar assumption on our terminal cost and terminal constraint set as in \cite{Amrit2011}. We use a straightforward extension of this assumption to periodic orbits. For systems with general operating regimes, similar ideas were followed in \cite{Dong2018}.
	\ass[Terminal conditions]{\label{ass:terminal_conditions}
		There exists a compact terminal region $\SX_\mathrm{f}$ satisfying $\varPi^\star_\SX(\SI_{[0,p^\star-1]})\subseteq\SX_\mathrm{f}\subseteq\SX$ such that for all $x\in\SX_\mathrm{f}$ there exists a control law $u_\mathrm{f}(x)\in\SU^1(x)$ such that $f(x,u_\mathrm{f}(x))\in\SX_\mathrm{f}$ and
		\begin{align}\label{eq:terminal_conditions}
			V_\mathrm{f}(x)\geq V_\mathrm{f}(f(x,u_\mathrm{f}(x))) + \ell(x,u_\mathrm{f}(x))-\ell_{\varPi^\star}.
		\end{align}
	}
	\vspace*{\belowdisplayskip}
	Under Assumption~\ref{ass:terminal_conditions}, standard MPC arguments show that the set inclusion $\SX_{N_0}\subseteq\SX_{N_1}$ holds for all $N_1\geq N_0$.
	In the literature, it is very often assumed that the optimal periodic orbit is minimal. To derive stability guarantees for non-minimal optimal periodic orbits, however, we impose an additional assumption on the terminal cost $V_\mathrm{f}$.
	\ass[Terminal cost on the optimal orbit]{\label{ass:terminal_cost_on_orbit}
		For all $(x,u)\in\varPi^\star(\SI_{[0,p^\star-1]})$, it holds that
		\begin{align}\label{eq:terminal_cost_on_orbit}
			V_\mathrm{f}(x)=\ell(x,u)-\ell_{\varPi^\star} + V_\mathrm{f}(f(x,u)).
		\end{align}
	}
	Note that Assumption~\ref{ass:terminal_cost_on_orbit} is trivially satisfied if $\varPi^\star$ is a steady state and, as we will show later in Proposition~\ref{prop:replace_assumption}, it is also satisfied if $\varPi^\star$ is minimal and Assumptions~\ref{ass:cont_comp}, \ref{ass:strict_dissipativity} and \ref{ass:terminal_conditions} hold. For non-minimal orbits, this implication is not true, as we illustrate in Example~\ref{ex:non-minimal_orbit} below.
	\example{\label{ex:non-minimal_orbit}
		Consider the system $x(k+1)=x(k)+u(k)$, with state and input constraints $\SX=[-1,1]$ and $\SU=[-1,1]$ and stage cost $\ell(x,u)=2(2\abs{x}-x^2-u^2)+1$. The system is optimally operated at the non-minimal periodic orbit $\varPi^\star(\SI_{[0,3]})=\{(-1,1),(0,1),(1,-1),(0,-1)\}$ with average cost $\ell_{\varPi^\star}=0$, as depicted in Figure~\ref{fig:non_minimal_orbit}. Indeed, for the storage function $\lb(x)=-x^2$, the rotated cost satisfies $\t\ell(x,u)=0$ for $(x,u)\in\varPi^\star(\SI_{[0,3]})$ and $\t\ell(x,u)>0$ otherwise, i.e., Assumption~\ref{ass:strict_dissipativity} holds. Now, let the terminal set and terminal cost be chosen as $\SX_\mathrm{f}\coloneq\varPi^\star_\SX(\SI_{[0,3]})$ and $V_\mathrm{f}(-1)\coloneq\frac{1}{2}$, $V_\mathrm{f}(0)\coloneq-\frac{1}{2}$ and $V_\mathrm{f}(1)\coloneq\frac{1}{2}$. For this choice, it is readily verified that Assumptions~\ref{ass:cont_comp}, \ref{ass:terminal_conditions} and \ref{ass:terminal_cost_on_orbit} are satisfied. However, for the same terminal set $\SX_\mathrm{f}$, let the terminal cost be chosen as $\h V_\mathrm{f}(-1)\coloneq1$, $\h V_\mathrm{f}(0)\coloneq V_\mathrm{f}(0)$ and $\h V_\mathrm{f}(1)\coloneq V_\mathrm{f}(1)$. Then again Assumption~\ref{ass:cont_comp} holds and also Assumption~\ref{ass:terminal_conditions} is satisfied with $u_\mathrm{f}(-1)\coloneq1$, $u_\mathrm{f}(0)\coloneq1$ and $u_\mathrm{f}(1)\coloneq-1$, but Assumption~\ref{ass:terminal_cost_on_orbit} is violated, since $-\frac{1}{2}=\h V_\mathrm{f}(0)<\ell(0,-1)+\h V_\mathrm{f}(f(0,-1))=0$.
		
	}
	
	
	\subsection{Candidate Lyapunov function} \label{sec:lyapunov_candidate}
	
	In order to construct a suitable candidate Lyapunov function, we define a rotated cost functional and a rotated optimal value function by replacing the stage cost $\ell$ by the rotated stage cost $\t\ell$ and by defining a suitable rotated terminal cost $\t V_\mathrm{f}$. For that we require the average storage and terminal cost over one period of the optimal periodic orbit $\varPi^\star$
	\[
	V_{\varPi^\star}\coloneq\frac{1}{p^\star}\sum_{i=0}^{p^\star-1}V_\mathrm{f}(\varPi_\SX^\star(i)),\quad \lb_{\varPi^\star}\coloneq\frac{1}{p^\star}\sum_{i=0}^{p^\star-1}\lb(\varPi_\SX^\star(i)).
	\]
	The rotated finite horizon cost functional then reads
	\begin{align*}
		\t J_N(x,u)\coloneq\sum_{k=0}^{N-1} \t\ell(x_u(k,x),u(k)) + \t V_\mathrm{f}(x_u(N,x)),
	\end{align*}
	with
	$
	\t V_\mathrm{f}(x):\SX_\mathrm{f}\to\R$ and $\t V_\mathrm{f}(x)\coloneq V_\mathrm{f}(x)+\lb(x)-V_{\varPi^\star}-\lb_{\varPi^\star}
	$
	and the corresponding rotated value function is defined by
	\[
	\t V_N(x)\coloneq\min_{\substack{u\in\SU^N(x) \\ x_u(N,x)\in\SX_\mathrm{f}}} \t J_N(x,u).
	\]
	Analogously to \cite[Lemma 14]{Amrit2011}, if Assumptions~\ref{ass:cont_comp}, \ref{ass:strict_dissipativity} and \ref{ass:terminal_conditions} are satisfied, it can be shown that the rotated and non-rotated optimization problems have the same set of minimizers, since
	$
	\t J_N(x,u)=J_N(x,u)-N\ell_{\varPi^\star}+\lb(x)-\lb_{\varPi^\star}-V_{\varPi^\star},
	$
	i,e., the objective functions $J_N(x,u)$ and $\t J_N(x,u)$ only differ by a constant. Consequently, it also holds that
	\begin{align}\label{eq:relation_accumulated_cost}
		\t V_N(x)=V_N(x)-N\ell_{\varPi^\star}+\lb(x)-\lb_{\varPi^\star}-V_{\varPi^\star}.
	\end{align}
	Next, we establish that $\t V_N$ is a suitable candidate Lyapunov function. In order to do so, we need several intermediate steps. First, we note that the inequality \eqref{eq:terminal_conditions} is equivalent to
	\begin{align} \label{eq:rot_term_cost_decrease}
		\t V_\mathrm{f}(x)\geq \t\ell(x,u_\mathrm{f}(x)) + \t V_\mathrm{f}(f(x,u_\mathrm{f}(x)))
	\end{align}
	for all $x\in\SX_\mathrm{f}$ and $u_\mathrm{f}(\cdot)$ as in Assumption~\ref{ass:terminal_conditions} and the equality in \eqref{eq:terminal_cost_on_orbit} is equivalent to
	\begin{align} \label{eq:rot_term_cost_on_orbit}
		\t V_\mathrm{f}(x)=\t \ell(x,u) + \t V_\mathrm{f}(f(x,u))
	\end{align}
	for all $(x,u)\in\varPi^\star(\SI_{[0,p^\star-1]})$. Both equivalences follow by simply plugging in the definitions of $\t V_\mathrm{f}$ and $\t\ell$. Secondly, we show that the rotated stage cost $\t\ell$ and the rotated terminal cost $\t V_\mathrm{f}$ vanish on the optimal periodic orbit.
	\lemma[Rotated costs on optimal orbit]{ \label{lem:rot_terminal_cost_on_orbit}
		Let Assumptions~\ref{ass:cont_comp}, \ref{ass:strict_dissipativity}, \ref{ass:terminal_conditions} and \ref{ass:terminal_cost_on_orbit} hold. Then, $\t\ell(x,u)=0$ and $\t V_\mathrm{f}(x)=0$ for all $(x,u)\in\varPi^\star(\SI_{[0,p^\star-1]})$.
	}
	\begin{proof}
		We first note that summing the rotated stage cost $\t\ell$ over one period of the optimal periodic orbit $\varPi^\star$ results in
		$
		\sum_{i=0}^{p^\star-1}\t\ell(\varPi^\star(i))=\sum_{i=0}^{p^\star-1}\left(\ell(\varPi^\star(i))-\ell_{\varPi^\star}\right)=0.
		$
		Hence, since $\t\ell(\cdot)$ is non-negative, necessarily $\t\ell(\varPi^\star(i))=0$ for all $i\in\SI_{[0,p^\star-1]}$. By Assumption~\ref{ass:terminal_cost_on_orbit} and \eqref{eq:rot_term_cost_on_orbit}, this implies that $\t V_\mathrm{f}(\varPi^\star_\SX(i))=\t V_\mathrm{f}(\varPi^\star_\SX(j))$ for all $i,j\in\SI_{[0,p^\star-1]}$, i.e., that $\t V_\mathrm{f}$ is constant on $\varPi^\star$. Since by definition $\sum_{i=0}^{p^\star-1}\t V_\mathrm{f}(\varPi^\star_\SX(i))=0$, it follows that $\t V_\mathrm{f}(x)=0$ for all $x\in\varPi^\star(\SI_{[0,p^\star-1]})$.
	\end{proof}
	Thirdly, we establish that the rotated terminal cost $\t V_\mathrm{f}$ is positive definite w.r.t. the optimal periodic orbit.
	\lemma[Rot. terminal cost is pos. def.]{\label{lem:rot_term_cost_pos_def}
		Let Assumptions~\ref{ass:cont_comp}, \ref{ass:strict_dissipativity}, \ref{ass:terminal_conditions} and \ref{ass:terminal_cost_on_orbit} be satisfied. There exist $\ubar\alpha,\obar\alpha\in\CK_\infty$ such that the rotated terminal cost $\t V_\mathrm{f}(x)$ satisfies for all $x\in\SX_\mathrm{f}$
		\[
		\ubar\alpha(\norm{(x,u_\mathrm{f}(x))}_{\varPi^\star})\leq\t V_\mathrm{f}(x)\leq\obar\alpha(\norm{(x,u_\mathrm{f}(x))}_{\varPi^\star}).
		\]
	}
	\vspace*{\belowdisplayskip}
	
	\begin{proof}
		Let $x\in\SX_\mathrm{f}$ and $u_\mathrm{f}(x)$ be the control law from Assumption \ref{ass:terminal_conditions} and abbreviate $u(k)\coloneq u_\mathrm{f}(x_{u_\mathrm{f}}(k,x))$ for $k\in\N_{\geq0}$. Then, by \eqref{eq:rot_term_cost_decrease},
		\begin{align}\label{eq:rotated_ineq_terminal}
			\begin{split}
				\t V_\mathrm{f}(x_{u}(k,x))\geq\t V_\mathrm{f}(x_{u}(k+1,x))+\t\ell(x_{u}(k,x),u(k)),
			\end{split}
		\end{align}
		for all $k\in\N_{\geq0}$. Since $\t\ell$ is non-negative, the sequence $\t V_\mathrm{f}(x_{u}(k,x))$ is non-increasing with $k$. Since $\t V_\mathrm{f}(x_{u}(k,x))$ is also bounded by compactness of $\SX_\mathrm{f}$ and continuity of $\t V_\mathrm{f}$, the sequence converges. Convergence of the sequence implies $\t\ell(x_{u}(k,x),u(k))\rightarrow0$ as $k\to\infty$ and consequently, by \eqref{eq:rotated_cost_pos}, also that $\norm{(x_{u}(k,x),u(k))}_{\varPi^\star}\rightarrow0$ as $k\rightarrow\infty$. Hence, by continuity of $\t V_\mathrm{f}$ and Lemma~\ref{lem:rot_terminal_cost_on_orbit}, we infer $\t V_\mathrm{f}(x_{u}(K,x))\rightarrow0$ as $K\to\infty$. Summing \eqref{eq:rotated_ineq_terminal} for $k\in\SI_{[0,K-1]}$, taking the limit $K\to\infty$ and using Assumption~\ref{ass:strict_dissipativity} then yields
		\begin{align*}
			\t V_\mathrm{f}(x)&\geq\limsup_{K\rightarrow\infty}\sum_{k=0}^{K-1}\t\ell(x_{u}(k,x),u(k)) + \t V_\mathrm{f}(x_{u}(K,x))\\
			&\geq\ubar\alpha_{\t\ell}(\norm{(x,u_\mathrm{f}(x))}_{\varPi^\star}).
		\end{align*}
		Therefore, we can choose $\ubar\alpha\coloneq\ubar\alpha_{\t\ell}\in\CK_\infty$ as lower bound. To find an upper bound $\t V_\mathrm{f}(x)\leq\obar\alpha(\norm{x}_{\varPi_\SX^\star})$, with $\obar\alpha\in\CK_\infty$, we can straightforwardly extend the proof of \cite[Lemma 4.3]{Khalil2002} to our setting, since $\t V_\mathrm{f}(x)$ is continuous on $x\in\SX_\mathrm{f}$ by Assumptions~\ref{ass:cont_comp} and \ref{ass:strict_dissipativity} and $\t V_\mathrm{f}(x)=0$ for $x\in\varPi_\SX^\star(\SI_{[0,p^\star-1]})$ by Lemma~\ref{lem:rot_terminal_cost_on_orbit}. Noting that $\norm{x}_{\varPi_\SX^\star}\leq\norm{(x,u)}_{\varPi^\star}$ for all $u\in\SU^1(x)$ concludes the proof.
	\end{proof}
	By Continuity of $\lb$ and compactness of $\SX$ there exists $\gamma_{\lb}\in\CK_\infty$ satisfying the bound $\abs{\lb(x)-\lb(\varPi^\star_\SX(i_x))}\leq\gamma_\lb(\norm{x}_{\varPi_\SX^\star})$ with $i_x$ from Definition~\ref{def:periodic_orbits}.  
	To establish $\t V_N$ as candidate Lyapunov function, a similar bound is required for $V_N$.
	\ass[bound on $V_N$]{\label{ass:bound_V_N}
		There exists $\gamma_V\in\CK_\infty$ such that for $i_x$ from Definition~\ref{def:periodic_orbits}, for each $N\in\N$ and each $x\in\SX_N$ it holds that
		$ 
		\abs{V_N(x)-V_N(\varPi^\star_\SX(i_x))}\leq\gamma_V(\norm{x}_{\varPi_\SX^\star}).
		$
	}
	If $\varPi^\star$ lies in the interior of $\SX_\mathrm{f}$, then there exists $\gamma_{\t V}\in\CK_\infty$ satisfying $|\t V_N(x)-\t V_N(\varPi^\star_\SX(i_x))|\leq\gamma_{\t V}(\norm{x}_{\varPi_\SX^\star})$, which can be proven by straightforwardly adapting the proofs of \cite[Proposition 2.17 and 2.18]{Rawlings2009}. Exploiting the relation between $V_N$ and $\t V_N$ then shows that also Assumption~\ref{ass:bound_V_N} holds. In case of terminal equality constraints, i.e., $\SX_\mathrm{f}=\varPi^\star_\SX(\SI_{[0,p^\star-1]})$, for example a local controllability condition such as \cite[Ass. 10]{Mueller2016} is sufficient for Assumption~\ref{ass:bound_V_N} to hold.
	We can now establish that $\t V_N$ is a suitable candidate Lyapunov function.
	\lemma[Candidate Lyapunov function]{\label{lem:Lyap_fcn}
		{Let Assumptions \ref{ass:cont_comp}, \ref{ass:strict_dissipativity}, \ref{ass:terminal_conditions}, \ref{ass:terminal_cost_on_orbit} and \ref{ass:bound_V_N} hold. There exist $\al_1,\al_2\in\CK_\infty$ such that the rotated value function $\t V_N(x)$ satisfies for all $x\in\SX_N$ 
		\vspace{-12pt}
		\begin{subequations}%
			\begin{align}%
				\al_1(\norm{x,\mu_N(x)}_{\varPi^\star})\leq \t V_N(x)\leq\al_2(\norm{x}_{\varPi_\SX^\star}), \label{eq:lyap_pos} \\
				\t V_N(f(x,\mu_N(x))) - \t V_N(x)\leq -\t\ell(x,\mu_N(x)).	\label{eq:lyap_neg}
			\end{align}
		\end{subequations}}
	}
	\vspace*{\belowdisplayskip}
	\begin{proof}
		The lower bound in \eqref{eq:lyap_pos} is satisfied by $\al_1\coloneq\ubar{\al}_{\t\ell}$, since $\t\ell(x,u)\geq\ubar{\al}_{\t\ell}(\norm{(x,u}_{\varPi^\star})\geq0$ for $u\in\SU^1(x)$ by \eqref{eq:rotated_cost_pos}, $\t V_\mathrm{f}(x)\geq0$ for $x\in\SX_\mathrm{f}$ by Lemma \ref{lem:rot_term_cost_pos_def} and hence
		\begin{align*}
			\t V_N(x)\geq\t\ell(x,u^\star_{N,x}(0))\stackrel{\eqref{eq:feedback_law},\eqref{eq:rotated_cost_pos}}{\geq}\ubar{\al}_{\t\ell}(\norm{(x,\mu_N(x))}_{\varPi^\star}).
		\end{align*}
		The upper bound in \eqref{eq:lyap_pos} is satisfied by $\alpha_2\coloneq\gamma_V+\gamma_\lb$, which follows by plugging \eqref{eq:relation_accumulated_cost} into the inequality in Assumption~\ref{ass:bound_V_N} and exploiting the bounds on $\lb$ from the paragraph above Assumption~\ref{ass:bound_V_N}. Finally, $\eqref{eq:lyap_neg}$ can be derived along the lines of the proof of \cite[Theorem 15]{Amrit2011}.
	\end{proof}
	We remark that \eqref{eq:lyap_neg} implies by \eqref{eq:relation_accumulated_cost} that also
	\begin{align}\label{eq:decrease}
		V_N(f(x,\mu_N(x)))-V_N(x)\leq-\ell(x,\mu_N(x))+\ell_{\varPi^\star}.
	\end{align}
	
	\subsection{Asymptotic Stability} \label{sec:asymptotic_stability}
	In the following we show that the optimal periodic orbit is asymptotically stable under the EMPC feedback law $\mu_N$. 
	\definition[Asymptotic stability of periodic orbits]{ \label{def:asy_stab}
		Let $\CS\subseteq\SX$ be a positively invariant set of system \eqref{eq:sys} under the control law $\mu:\CS\rightarrow\SU$. A $p$-periodic orbit $\varPi$ is called locally stable in $\CS$ under the feedback $\mu$, if for all $\eps>0$ there exists $\de>0$ such that $\norm{x}_{\varPi_\SX}<\delta$ and $x\in\CS$ imply $\norm{(x_{\mu}(k,x),\mu(x_{\mu}(k,x)))}_\varPi<\eps$ for all $k\in\N_{\geq0}$. It is called asymptotically stable in $\CS$ with region of attraction $\CS$, if additionally $\norm{(x_{\mu}(k,x),\mu(x_{\mu}(k,x)))}_\varPi\rightarrow0,$ as $k\to\infty$ for all $x\in\CS$.
	}
	In contrast to asymptotic stability of a steady state, in Definition~\ref{def:asy_stab}, we not only impose conditions on the state but also on the feedback law. This ensures that we not only have some form of set stability but a stronger version of stability, i.e., that the state trajectory is guaranteed to move along the possibly non-minimal periodic orbit, which does however not imply that the resulting state trajectory is periodic, as explained in the following example.
	\example{\label{ex:non-periodic}
		Consider the non-minimal optimal periodic orbit $\varPi^\star$ from Example~\ref{ex:non-minimal_orbit} (cf. also Fig.~\ref{fig:non_minimal_orbit}). The optimal orbit $\varPi^\star$ can be decomposed into two sub-orbits $\varPi^1(\SI_{[0,1]})=\{(0,-1),(-1,1)\}$ and $\varPi^2(\SI_{[0,1]})=\{(0,1),(1,-1)\}$. Then, when starting on the optimal periodic orbit $\varPi^\star$, asymptotic stability as in Definition~\ref{def:asy_stab} guarantees that the state–input trajectory can be expressed as a concatenation of $\varPi^1$ and $\varPi^2$, but this concatenation need not be periodic.
	}
	\theorem[Asy. Stability]{ \label{thm:asy_stab}
		Let Assumptions \ref{ass:cont_comp}, \ref{ass:strict_dissipativity}, \ref{ass:terminal_conditions}, \ref{ass:terminal_cost_on_orbit} and \ref{ass:bound_V_N} hold and let system \eqref{eq:sys} be controlled by the optimal feedback law \eqref{eq:feedback_law}. Then, $\varPi^\star$ is asymptotically stable according to Definition~\ref{def:asy_stab} with region of attraction $\SX_N$.
	}
	\begin{proof} The following arguments are adapted from \cite[Theorem B.11]{Rawlings2009}.
		\emph{Stability:} Consider $\al_1,\al_2\in\CK_\infty$ from Lemma \ref{lem:Lyap_fcn}. For any $\eps>0$, let $\de\coloneq\al^{-1}_2(\al_1(\eps))$. Suppose $x\in\SX_N$ and $\norm{x}_{\varPi_\SX^\star}<\de$. Then, \eqref{eq:lyap_pos} implies $\t V_N(x)\leq\al_1(\eps)$. By \eqref{eq:lyap_neg} we infer that $\t V_N(x_{\mu_N}(k,x))\leq \t V_N(x)$  for $k\in\N$ and hence by \eqref{eq:lyap_pos} that $\norm{(x_{\mu_N(k,x)},\mu(x_{\mu_N(k,x)}))}_{\varPi^\star}\leq\al_1^{-1}(\t V_N(x))\leq\eps$ for $k\in\N$. 
		
		\emph{Convergence:} Note that for all $x\in\SX_N$, the sequence $(\t V_N(x_{\mu_N}(k,x)))_{k\in\N}$ is non-increasing and bounded from below. Hence, both $\t V_N(x_{\mu_N}(k+1,x))$ and $\t V_N(x_{\mu_N}(k,x))$ converge to the same limit, i.e., $\t V_N(x_{\mu_N}(k+1,x))-\t V_N(x_{\mu_N}(k,x))\to0$ as $k\to\infty$. With  ${x_{\mu_N}(k+1,x)}=f(x_{\mu_N}(k,x),\mu_N(x_{\mu_N}(k,x)))$, we infer by \eqref{eq:lyap_neg} that $\t\ell(x_{\mu_N}(k,x),\mu_N(x_{\mu_N}(k,x)))\to0$ as $k\to\infty$. By \eqref{eq:rotated_cost_pos}, this implies\! $\|(x_{\mu_N}(k,x),\mu_N(x_{\mu_N}(k,x)))\|_{\varPi^\star}\!\to0$ as $k\to\infty$.
	\end{proof}
	Theorem~\ref{thm:asy_stab} naturally extendeds the well-known stability result for EMPC for systems which are optimally operated at steady states \cite{Amrit2011}, to systems with periodic optimal operation. Compared to \cite{Dong2018}, we use a slightly stronger dissipativity notion, since the distance to the optimal periodic orbit in \eqref{eq:rotated_cost_pos} depends on the pair $(x,u)$ and not only on $x$. However, by doing this, we do not need to assume continuity of the value function $V_N(x)$ and the feedback law $\mu_N(x)$ in a neighborhood of the orbit, which excludes general non-minimal optimal periodic orbits. In contrast, our stability result does not require minimality of the optimal periodic orbit, which is often the case for other EMPC schemes dealing with periodic orbits \cite{Mueller2016}, \cite{Schwenkel2024}. Instead, we need Assumption~\ref{ass:terminal_cost_on_orbit}. However, if the optimal periodic orbit is in fact minimal, then Assumption~\ref{ass:terminal_cost_on_orbit} is implied by Assumptions~\ref{ass:cont_comp}, \ref{ass:strict_dissipativity} and \ref{ass:terminal_conditions} and therefore it is in this case not restrictive.
	\proposition[Relaxation for minimal orbits]{\label{prop:replace_assumption}
		Let Assumptions~\ref{ass:cont_comp}, \ref{ass:strict_dissipativity} and \ref{ass:terminal_conditions} hold. Then, Assumption~\ref{ass:terminal_cost_on_orbit} holds if the optimal periodic orbit $\varPi^\star$ is minimal.
	}
	\begin{proof}
		Let $\t V_\mathrm{min}\coloneq\min_{x\in\SX_\mathrm{f}}\t V_\mathrm{f}(x)$ and $\SX_\mathrm{min}\coloneq\{x\in\SX_\mathrm{f}\,|\,\t V(x)=\t V_\mathrm{min}\}$. Let $x\in\SX_\mathrm{min}$ and $u_\mathrm{f}(\cdot)$ be from Assumption~\ref{ass:terminal_conditions}. Then, we infer by \eqref{eq:rot_term_cost_decrease} that
		\begin{align*}
			\t\ell(x,u_\mathrm{f}(x))+\t V_\mathrm{f}(f(x,u_\mathrm{f}(x)))\leq\t V_\mathrm{f}(x)=\t V_\mathrm{min}.
		\end{align*}
		Hence, since $\t V_\mathrm{f}(f(x,u_\mathrm{f}(x)))\geq\t V_\mathrm{min}$ by optimality, and $\t \ell(\cdot)$ is non-negative, necessarily $\t V_\mathrm{f}(f(x,u_\mathrm{f}(x)))=\t V_\mathrm{min}$ and $\t\ell(x,u_\mathrm{f}(x))=0$, which implies $x\in\varPi^\star_\SX(\SI_{[0,p^\star-1]})$ by \eqref{eq:rotated_cost_pos}. By minimality of $\varPi^\star$, there exists a unique $i\in\SI_{[0,p^\star-1]}$ such that $(x,u_\mathrm{f}(x))=\varPi^\star(i)$. Hence, by the arguments above, $\t V_\mathrm{f}(\varPi^\star_\SX(i))=\t V_\mathrm{f}(\varPi^\star_\SX([i+1]_{p^\star}))=\t V_\mathrm{min}$. By induction it follows that $\t V_\mathrm{f}(x)=V_\mathrm{min}$ for all $x\in\varPi^\star_\SX$, which proves \eqref{eq:rot_term_cost_on_orbit} and therefore also that Assumption~\ref{ass:terminal_cost_on_orbit} holds.
	\end{proof}
	
	\subsection{Numerical Example}
	
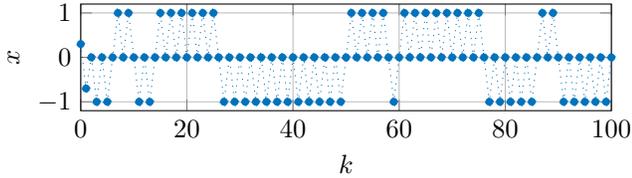
\begin{figure}[tp]
		\medskip
		\begin{tikzpicture}
			\begin{axis} [grid=major, legend pos=south east, legend style={font=\tiny,  legend columns=-1},
				ymin=-1.2, ymax=1.2, ylabel={$x$},
				xmin=0, xmax=100, xlabel={$k$}, x label style={at={(axis description cs:0.5,-0.2)}}, 
				height=3.0cm, width = \columnwidth,
				cycle list name=convPlot]
				\addplot table [x index=0, y index=1, col sep=comma] {figs/state_trajectory_ic2.txt};
			\end{axis}
		\end{tikzpicture}
		\caption{Closed-loop state trajectory $x_{\mu_{5}}(k,0.3)$ over the time index $k$.}
		\label{fig:state-trajectory}
	\end{figure}

	In this subsection, we illustrate the stability results of subsection~\ref{sec:asymptotic_stability} with a numerical example. Consider again the system and state-input constraints from Example~\ref{ex:non-minimal_orbit} and the corresponding non-minimal optimal periodic orbit $\varPi^\star$. Let also the terminal set and terminal cost be given by $\SX_\mathrm{f}$ and $V_\mathrm{f}$ from Example~\ref{ex:non-minimal_orbit}. We employ the EMPC scheme from Secion~\ref{sec:MPC-scheme} with prediction horizon $N=5$. The state trajectory resulting from the EMPC feedback law $\mu_{5}(x)$ is depicted in Figure~\ref{fig:state-trajectory} for $x(0)=0.3$. The controller steers the system's state to $\varPi^\star$ in two steps and then remains on $\varPi^\star$, however, as explained in Example~\ref{ex:non-periodic}, the state-trajectory is not necessarily periodic. Note that the optimal feedback law $\mu_5(x)$ is not unique at $x=0$, and hence it is also not continuous at this point, which implies that the result in \cite{Dong2018} cannot be applied here, whereas Theorem~\ref{thm:asy_stab} is applicable\footnote{Note that Assumption~\ref{ass:bound_V_N} holds, since the system is locally controllable (cf. \cite[Assumption 10]{Mueller2016})}.
	\pgfplotscreateplotcyclelist{convPlot}{dotted,RoyalBlue,line width=0.4pt,mark size=1.5pt,mark=*\\}
	
	\section{Performance Guarantees} \label{sec:performance}
	
	In this section, we extend the non-averaged performance guarantees in \cite{Gruene2015} to systems with periodic optimal operation. The advantage of non-averaged performance compared to averaged performance is that transient costs do not vanish asymptotically. From now on we assume w.l.o.g. that $\ell_{\varPi^\star}=0$, since the stage cost can always be shifted accordingly without changing the minimizer. For systems with optimal operation at steady state, the closed-loop cost over a horizon $K\in\N$ is defined as the sum over the stage costs in closed-loop operation, i.e., the closed-loop cost is defined as $J_K^1(x,\mu)$, with
	\begin{align*}
		J^1_K(x,u)\coloneq\sum_{k=0}^{K-1}\ell(x_u(k,x),u(k)).
	\end{align*}
	However, for systems with periodic optimal operation, this sum does in general not converge for $K\to\infty$. Thus, we use Cesàro summation to define optimal costs for systems with periodic optimal operation, since it is well-known that Cesàro sums converge for periodic sequences with zero mean~\cite{Hardy1949}. Accordingly to \cite[Definition 4]{Mair25}, we define the Cesàro cost
	\begin{align*}
		J^\mathrm{ces}_K(x,u)&\coloneq\frac{1}{K}\sum_{n=0}^{K-1} J^1_{n+1}(x,u) \\
		&=\sum_{k=0}^{K-1}\left(1-\frac{k}{K}\right)\ell(x_u(k,x),u(k))
	\end{align*}	
	for the finite horizon and for the infinite horizon $J^\mathrm{ces}_\infty(x,u)\coloneq\limsup_{K\rightarrow\infty}J^\mathrm{ces}_K(x,u)$. The closed-loop cost under the feedback law $\mu_N$ is then defined as $J^\mathrm{cl}_K(x,\mu_N)\coloneq J^\mathrm{ces}_K(x,\mu_N)$ and $J^\mathrm{cl}_\infty(x,\mu_N)\coloneq J^\mathrm{ces}_\infty(x,\mu_N)$. For a fixed control law $\mu(x)$, the stage cost sequence can be regarded as a sequence $\ell_k(x)\coloneq\ell(x_\mu(k,x),\mu(x_\mu(k,x)))$ for $k\in\N$ that only depends on $x$. Therefore, well-known results regarding the relation of limits of summation methods directly yield that if the limit $J^1_\infty(x,\mu)=\lim_{N\rightarrow\infty}J^1_N(x,\mu)$ exists, then $J^1_\infty(x,\mu)=J^\mathrm{ces}_\infty(x,\mu)=J^\mathrm{cl}_\infty(x,\mu)$ \cite{Hardy1949}. Extending the definitions of the closed-loop cost via the Cesàro summation for systems with periodic optimal operation is therefore in that sense consistent with the definitions in \cite{Gruene2015}. Analogously to \cite{Gruene2015}, we furthermore define an unconditioned optimal cost without terminal constraint and terminal cost. Again we replace ordinary summation by a Cesàro sum and therefore
	\[
	V_N^\mathrm{uc}(x)\coloneq\min_{u\in\SU^N(x)} J^\mathrm{ces}_N(x,u)
	\]
	for the finite horizon and for the infinite horizon
	\[
	V_\infty^\mathrm{uc}(x)\coloneq\min_{u\in\SU^\infty(x)} \limsup_{K\rightarrow\infty}J^\mathrm{ces}_K(x,u).
	\]
	The corresponding rotated quantities $\t V^\mathrm{uc}_N$ and $\t J^\mathrm{ces}_N$ are defined using $\t\ell$ instead of $\ell$. Note that, by \cite[Lemma 21]{Schwenkel2024}, 
	\begin{align} \label{eq:relation_rotated_nonrotated}
		\begin{split}
			\t J^\mathrm{ces}_N(x,u)=J^\mathrm{ces}_N(x,u)
			+\lb(x) 
			-&\frac{1}{N}\sum_{k=1}^{N}\lb(x_u(k,x))
		\end{split}
	\end{align}
	describes the relation of the rotated and original finite horizon cost for $\ell_{\varPi^\star}=0$.
	We stress that if the unconditioned infinite horizon cost defined by ordinary summation as in \cite{Gruene2015} is well-defined in the sense that for all $x\in\SX_\infty$ it holds that 
	\begin{align*}
		V^1_\infty(x)&\coloneq\inf_{u\in\SU^\infty(x)} \limsup_{N\rightarrow\infty} J^1_N(x,u) \\
		&=\inf_{u\in\SU^\infty(x)} \liminf_{N\rightarrow\infty} J^1_N(x,u),
	\end{align*}
	then, under a suitable controllability and reachability assumption, also $V^1_\infty(x)=V^\mathrm{uc}_\infty(x)$ holds \cite[Theorem 3]{Mair25}.

	\subsection{Preliminary results} \label{sec:preliminary_resuluts}
	
	In this subsection we derive some preliminary results which are required to prove non-averaged performance guarantees for the EMPC scheme.
	\ass[Minimality of $\varPi^\star$]{\label{ass:minimality}The optimal periodic orbit $\varPi^\star$ is minimal.}
	We stress that minimality of $\varPi^\star$ was not needed for the previous asymptotic stability result, but it is crucial for deriving the non-averaged performance bound. We can now derive bounds on the infinite-horizon closed-loop cost $J^\mathrm{cl}_\infty$ and the unconditioned cost $V^\mathrm{uc}_\infty$.
	\lemma[Upper bound on closed-loop cost]{\label{lem:Upper-bound_on_closed-loop_performance}
		Let Assumptions \ref{ass:cont_comp}, \ref{ass:strict_dissipativity}, \ref{ass:terminal_conditions}, \ref{ass:bound_V_N} and \ref{ass:minimality} hold and let $\ell_{\varPi^\star}=0$. Then, $J^\mathrm{cl}_\infty(x,\mu_N)\leq V_N(x)-V_{\varPi^\star}$ for all $N\in\N$ and all $x\in\SX_N$.
	}
	\begin{proof}
		Using \eqref{eq:decrease} and $\ell_{\varPi^\star}=0$ yields for any $K>0$
		\begin{align*}
			J^\mathrm{cl}_K(x,u_x)&=\sum_{k=0}^{K-1} (1-\frac{k}{K})\ell(x_{\mu_N}(k,x),\mu_N(x_{\mu_n}(k,x))) \\
			&\leq V_N(x)-\frac{1}{K}\sum_{k=1}^{K} V_N(x_{\mu_N}(k,x)).
		\end{align*}
		Since we have asymptotic stability, we can now follow the proof of \cite[Lemma 15.4]{Schwenkel2024a} with $\beta(x)\coloneq1-x$, to infer the existence of $\delta\in\CL$ such that for all $x\in\SX_N$
		\[
		\abs{\frac{1}{p^\star}\sum_{i=0}^{p^\star-1}V_N(\varPi^\star_\SX(i))-\frac{1}{K}\sum_{k=1}^{K} V_N(x_{\mu_N}(k,x))}\leq\delta(K)
		\]
		holds\footnote{The proof of \cite[Lemma 15.4]{Schwenkel2024a} exploits a weak turnpike property of the state-input trajectory. We did not show such a property here. However, asymptotic stability is an even stronger property, since it implies that for all $\eps>0$ there exists $N_0\in\N$ such that the state-input trajectory stays inside an $\eps$-neighborhood around the optimal periodic orbit for all $N\geq N_0$.}. It is easy to verify that $V_{\varPi^\star}=\frac{1}{p^\star}\sum_{i=0}^{p^\star-1}V_N(\varPi^\star_\SX(i))$ and therefore $\limsup_{K\rightarrow\infty}J^\mathrm{cl}_K(x,\mu_N)\leq V_N(x)-V_{\varPi^\star}$.
	\end{proof}
	\lemma[Bounds on $V^\mathrm{uc}_\infty$]{ \label{lem:bounds_on_Vuc}
		Let Assumptions \ref{ass:cont_comp}, \ref{ass:strict_dissipativity}, \ref{ass:terminal_conditions}, \ref{ass:bound_V_N} and \ref{ass:minimality} hold and let $\ell_{\varPi^\star}=0$. Then there exist $B,D>0$ such that for all $x\in\CX_\infty$ it holds that
		\[
		-D\leq V^\mathrm{uc}_\infty(x)\leq\gamma_V(\norm{x}_{\varPi_\SX^\star})+B-	V_{\varPi^\star}.
		\]
	}
	\vspace*{\belowdisplayskip}
	\begin{proof}
		For any $N\in\N$, let $x\in\SX_N$ and note that $J^\mathrm{ces}_\infty(x,\mu_N)=J^\mathrm{cl}_\infty(x,\mu_N)$. Hence, by optimality, Lemma~\ref{lem:Upper-bound_on_closed-loop_performance} and Assumption~\ref{ass:bound_V_N}, it holds that
		\begin{align*}
			V^\mathrm{uc}_\infty(x)&\leq J^\mathrm{ces}_\infty(x,\mu_N) \leq V_N(x)-V_{\varPi^\star}\\ 
			&	\leq\gamma_V(\norm{x}_{\varPi^\star_\SX})+V_N(\varPi^\star_\SX(i_x))-V_{\varPi^\star}.
		\end{align*}
		Now construct a candidate input $\bar u(k)\coloneq(\varPi_\SU^\star([k+i_x]_{p^\star}))$ for $k\in\SI_{[0,N-1]}$, which follows the optimal orbit $\varPi^\star$. Then, by optimality and since the average cost on $\varPi^\star$ satisfies $\ell_\varPi^\star=0$, we infer for any $N$
		\begin{align*}
			V_N(\varPi^\star_\SX&(i_x))\leq J_N(\varPi^\star_\SX(i_x),\bar u) \\
			\leq& \;V_\mathrm{f}(\varPi^\star_\SX([i_x+N]_{p^\star}))+(p^\star-1)\max_{(x,u)\in\SX\times\SU}\ell(x,u).
		\end{align*}
		Hence, by Assumption~\ref{ass:cont_comp}, there exists $B<\infty$ such that $V_N(\varPi^\star_\SX(i_x))\leq B$.
		For the lower bound, note that $\t J^\mathrm{ces}_N(x,u)\geq0$. By \eqref{eq:relation_rotated_nonrotated} and boundedness of $\lb$, there exists $D>0$ such that $J^\mathrm{ces}_N(x,u)\geq-D$ for all $x$, $u$ and $N$. This implies $V^\mathrm{uc}_\infty(x)\geq-D$ Since $N$ was arbitrary, both bounds also hold for $N\rightarrow\infty$ and hence for $x\in\CX_\infty$.
	\end{proof}
	As it is shown in \cite{Schwenkel2024}, optimal trajectories resulting from minimizing the Cesàro cost $J^\mathrm{ces}_K$ satisfy a weak turnpike property w.r.t. to the optimal periodic orbit $\varPi^\star$. It is called a weak turnpike, because infinitely many points of the trajectory may lie outside of an $\eps$-neighborhood of $\varPi^\star$. However, the ratio of points outside the $\eps$-neighborhood and the total amount of points is still zero asymptotically.
	\lemma[Weak Turnpike property]{ \label{lem:weak_Turnpike}
		Let Assumptions~\ref{ass:cont_comp} and \ref{ass:strict_dissipativity} hold. Then there exists $C>0$ such that for each $x\in\SX$, $\delta>0$, $K\in\N$, each control sequence $u\in\SU^K(x)$ satisfying $J^\mathrm{ces}_K(x,u)\leq\tfrac{K+1}{2}\ell_{\varPi^\star}+\delta$ and each $\eps>0$ the value $Q_\eps(x)\coloneq\#\{k\in\SI_{[0,K-1]}\,|\,\norm{x_u(k,x),u(k)}_{\varPi^\star}\leq\eps\}$ satisfies the inequality $Q_\eps(x)\geq K-\tfrac{\sqrt{2(\de+C)K}}{\ubar\alpha_{\t\ell}(\eps)}$.
	}
	\begin{proof}
		From $u\in\SU^K(x)$ satisfying $J^\mathrm{ces}_K(x,u)\leq\tfrac{K+1}{2}\ell_{\varPi^\star}+\delta$, boundedness of $\lb$ and \eqref{eq:relation_rotated_nonrotated} we infer there exists $C>0$ such that $\t J^\mathrm{ces}_K(x,u)\leq C+\de$. The rest follows from \cite[Lemma 23]{Schwenkel2024}.
	\end{proof}
	
	
	While the weak turnpike property established in Lemma~\ref{lem:weak_Turnpike} characterizes the behavior of optimal trajectories on a finite horizon, we can also characterize the behavior of almost optimal trajectories on an infinite horizon.
	
	\lemma[Behavior of close to optimal trajectories]{\label{lem:clost_to_optimal_trajectories}
		Let Assumptions~\ref{ass:cont_comp}, \ref{ass:strict_dissipativity}, \ref{ass:terminal_conditions}, \ref{ass:bound_V_N} and \ref{ass:minimality} hold and let $\ell_{\varPi^\star}=0$. Then there exists $\sg\in\CL$ such that for any $x\in\SX_\infty$, any $u\in\SU^\infty(x)$ with $J^\mathrm{ces}_\infty(x,u)\leq V^\mathrm{uc}_\infty(x)+1$ and any $K\in\N$ and $q\in\N$ there is $k\in\N$ with $q\leq k\leq K+q$ such that $\norm{(x_u(k,x),u(k))}_{\varPi^\star}\leq\sg(K)$.
	}
	\begin{proof}
		We start with the case $q=0$. By Assumption~\ref{ass:cont_comp} and Lemma \ref{lem:bounds_on_Vuc}, we infer that there exists $M<\infty$ such that $V^\mathrm{uc}_\infty(x)\leq M$ and hence, by assumption, $J^\mathrm{ces}_\infty(x,u)\leq M+1$. From \eqref{eq:relation_rotated_nonrotated}, boundedness of $\lb$, the fact that $\ell_\varPi^\star=0$, and Lemma~\ref{lem:monotonic_increase_Juc} in the Appendix, we then infer that there exists $C_1>0$ such that $\t J^\mathrm{ces}_K(x,u)\leq\t J^\mathrm{ces}_\infty(x,u)\leq C_1+M+1$ for any $K\in\N$. Again exploiting \eqref{eq:relation_rotated_nonrotated} and boundedness of $\lb$ then implies that there exists $C_2>0$ such that $ J^\mathrm{ces}_K(x,u)\leq C_1+C_2+M+1$ for any $K\in\N$. Hence, we can apply Lemma~\ref{lem:weak_Turnpike} with $\delta\coloneq C_1+C_2+M+1$. Without loss of generality, we can assume $C_1+C_2=C$. Then, choosing $\eps=\sg(K)\coloneq\ubar\alpha_{\t\ell}^{-1}(\tfrac{\sqrt{2(2C+M+1)K}}{K-1})$ in Lemma~\ref{lem:weak_Turnpike} leads to $Q_\eps\geq1$, which concludes the proof. For arbitrary $q\in\N$, we can use Lemma~\ref{lem:bound_propagation} in the Appendix to infer that we can just replace $x$ by $x_u(q,x)$ in the arguments above.
	\end{proof}
	Furthermore, the lower bound on $V^\mathrm{uc}_\infty$ from Lemma~\ref{lem:bounds_on_Vuc} can be improved.
	\lemma[Lower bound $V^\mathrm{uc}_\infty$]{\label{lem:lower_bound_Vuc}
		Let Assumptions~\ref{ass:cont_comp}, \ref{ass:strict_dissipativity}, \ref{ass:terminal_conditions}, \ref{ass:bound_V_N} and \ref{ass:minimality} hold and let $\ell_{\varPi^\star}=0.$ Then, $V^\mathrm{uc}_\infty(x)\geq-\lb(x)+\lb_{\varPi^\star}$ for all $x\in\SX_\infty$.
	}
	\begin{proof}
		Consider an arbitrary $\eps\in(0,1)$ and let $u\in\SU^\infty(x)$ such that $J^\mathrm{ces}_\infty(x,u)\leq V^\mathrm{uc}_\infty(x)+\eps$. Then, \eqref{eq:relation_rotated_nonrotated} and the definition of $V^\mathrm{uc}_\infty$ imply
		\[
		V^\mathrm{uc}_\infty(x)+\eps\geq\limsup_{K\rightarrow\infty} \t J^\mathrm{ces}_K(x,u)-\lb(x)+\frac{1}{K}\sum_{k=1}^{K}\lb(x_u(k,x)).
		\]
		Since $V^\mathrm{uc}_\infty(x)$ is upper bounded by Lemma~\ref{lem:bounds_on_Vuc} and compactness of $\SX$, we can apply Lemma~\ref{lem:weak_Turnpike} to infer that a weak turnpike property holds for the chosen state-input trajectory. Hence, we infer by \cite[Lemma 15.4]{Schwenkel2024a} $\lim_{K\rightarrow\infty}\tfrac{1}{K}\sum_{k=1}^{K}\lb(x_u(k,x))=\lb_{\varPi^\star}$ and therefore, since $\eps$ was arbitrary, $V^\mathrm{uc}_\infty(x)\geq-\lb(x)+\lb_{\varPi^\star}$.
	\end{proof}
	
	\subsection{Infinite-horizon non-averaged performance} \label{sec:non-averaged}
	
	We now use the previous section's results to prove a non-averaged performance bound on the closed-loop cost. To avoid additional technicalities, we assume that for a large enough horizon, a neighborhood of the optimal periodic orbit is contained in the feasible set of the optimization problem~\eqref{eq:finite_horizon_value_function}. This assumption was also imposed in \cite{Gruene2015} for the same reason.
	
	\ass[Feasibility close to $\varPi^\star$]{\label{ass:feasibility_close_to_orbit}
		There exists $N_\eta\in\N$ and $\eta>0$ such that $\norm{x}_{\varPi^\star_\SX}\leq\eta$ implies $x\in\SX_{N_\eta}$.
	}
	
	Using the definition of non-averaged closed-loop costs via the Cesàro summation, we now extend the non-averaged performance result in \cite{Gruene2015} for systems with optimal operation at steady state to systems with optimal periodic operation. The following theorem shows that the closed-loop cost approximates the unconditioned infinite horizon cost up to an error vanishing with increasing horizon length $N$. 
	
	\theorem[Non-averaged performance]{\label{thm:non-averaged_performance}
		Let Assumptions~\ref{ass:cont_comp}, \ref{ass:strict_dissipativity}, \ref{ass:terminal_conditions}, \ref{ass:bound_V_N}, \ref{ass:minimality} and \ref{ass:feasibility_close_to_orbit} hold and let $\ell_{\varPi^\star}=0$. Then there exists $\delta\in\CL$ such that for all $N\in\N$ and $x\in\SX_N$
		\[
		J^\mathrm{cl}_\infty(x,\mu_N)\leq V_N(x)-V_{\varPi^\star}\leq V^\mathrm{uc}_\infty(x)+\de(N).
		\]
	}
	\vspace*{\belowdisplayskip}
	
	\begin{proof}
		The first inequality is given by Lemma~\ref{lem:Upper-bound_on_closed-loop_performance}.\\
		For proving the second inequality, note that for small $N$, the value functions  $V_N$ and $V^\mathrm{uc}_\infty$ and are bounded and hence the second inequality can always be satisfied by choosing $\de(N)$ sufficiently large without violating the condition $\de\in\CL$. Hence, we prove the second inequality only for sufficiently large $N$. Pick $\eps\in(0,1)$ and choose an admissible control $u_\eps\in\SU^\infty(x)$ such that $J^\mathrm{ces}_\infty(x,u_\eps)\leq V^\mathrm{uc}_\infty(x)+\eps$. Then, for $N\in\N$, applying Lemma~\ref{lem:clost_to_optimal_trajectories} with $K=\floor{\tfrac{N}{2}}$ implies the existence of $k\in\SI_{[0,K]}$ and $\sg\in\CL$ satisfying $\norm{(x_{u_\eps}(k,x),u_\eps(k))}\leq\sg(K)$. Therefore, for $K$, i.e., $N$, large enough such that $\sg(K)<\eta$ and $K\geq N_\eta$ with $\eta$ and $N_\eta$ from Assumption~\ref{ass:feasibility_close_to_orbit}, we have $x_{u_\eps(k,x)}\in\SX_{N_\eta}\subseteq\SX_K\subseteq \SX_{N-k}$, since $N-k\geq K\geq N_\eta$. Next, by Lemma~\ref{lem:separation_infinite-horizon_costs} in the Appendix, optimality and Lemma~\ref{lem:lower_bound_Vuc}, we have
		\begin{align*}
			V^\mathrm{uc}_\infty(x)+\eps&\geq J^1_k(x,u_\eps) + J^\mathrm{ces}_\infty(x_{u_\eps}(k,x),u_\eps(k+\cdot)) \\
			&\geq J^1_k(x,u_\eps)+V^\mathrm{uc}_\infty(x_{u_\eps}(k,x)) \\
			&\geq J^1_k(x,u_\eps) - \lb(x_{u_\eps}(k,x)) + \lb_{\varPi^\star}.
		\end{align*}
		Note that by relation \eqref{eq:relation_accumulated_cost} and the inequalities in \eqref{eq:lyap_pos}, $V_N(\varPi_\SX^\star(i_x))+\lb(\varPi_\SX^\star(i_x))-\lb_{\varPi^\star}=V_{\varPi^\star}$ for all $x\in\SX_N$.
		Hence, by exploiting optimality, the inequality above and the bounds on $V_N$ and $\lb$ in terms of $\gamma_V$ and $\gamma_\lb$ from Assumption~\ref{ass:bound_V_N} and the paragraph above, we infer
		\begin{align*}
			V_N(x)\leq J^1_k&(x,u_\eps)+V_{N-k}(x_{u_\eps}(k,x)) \\
			\leq V^\mathrm{uc}_\infty(x) &+\eps + \lb(x_{u_\eps}(k,x))+ V_{N-k}(x_{u_\eps}(k,x))-\lb_{\varPi^\star} \\
			\leq V^\mathrm{uc}_\infty(x) &+\eps + \gamma_\lb(\sg(K)) + \gamma_V(\sg(K)) + V_{\varPi^\star}.
		\end{align*}
		The assertion $\de(N)\coloneq\gamma_\lb(\sg(\floor{\tfrac{N}{2}})) + \gamma_V(\sg(\floor{\tfrac{N}{2}}))$ concludes the proof, since $\eps>0$ was arbitrary.
	\end{proof}

	\section{Cesàro Cost as Terminal Cost} \label{sec_cesaro_terminal}
	
	In the previous section, we derived a non-averaged performance bound in terms of Cesàro costs of an EMPC scheme with terminal conditions. This motivates to also examine the suitability of Cesàro costs $V^\mathrm{ces}_N(x)\coloneq V^\mathrm{uc}_N(x)$ as a terminal cost for systems with periodic optimal operation, i.e. we assume that Assumption~\ref{ass:strict_dissipativity} is satisfied. We first note that if $f$ and $\ell$ are continuous, then for any compact subset $\SX_\mathrm{c}\subseteq\SX_\infty$ and for all $x\in\SX_\mathrm{c}$, it holds that 
	\begin{equation}\label{eq:cesaro_iteration}
		V^\mathrm{ces}_N(x)= \min_{u\in\SU^1(x)} \ell(x,u) + \frac{N-1}{N} V^\mathrm{ces}_{N-1}(f(x,u))
	\end{equation} 
	for all $N\in\N$ \cite[Theorem 8]{Mair25}. Therefore, if the limit $V^\mathrm{ces}_\infty(x)\coloneq\lim_{N\rightarrow\infty} V^\mathrm{ces}_N(x)$ exists for all $x\in\SX_\mathrm{c}$, which directly implies $\ell_{\varPi^\star}=0$, then there exists $u\in\SU^1(x)$ such that
	$
	V^\mathrm{ces}_\infty(x) = \ell(x,u) + V^\mathrm{ces}_\infty(f(x,u))
	$
	for all $x\in\SX_\mathrm{c}$. For $\SX_\mathrm{c}=\SX_\mathrm{f}$ we hence conclude that $V^\mathrm{ces}_\infty$ satisfies Assumption~\ref{ass:terminal_conditions}. 
	Furthermore, it can be shown that $V^\mathrm{ces}_\infty$ also satisfies Assumption~\ref{ass:terminal_cost_on_orbit}, as formalized in Proposition~\ref{prop:cesaro_ass5} below.
	\proposition[Cesàro cost satisfies Ass.~\ref{ass:terminal_cost_on_orbit}]{ \label{prop:cesaro_ass5}
		Let $\SX_\mathrm{f}$ be compact and satisfy $\varPi_\SX^\star(\SI_{[0,p^\star-1]})\subseteq\SX_\mathrm{f}\subseteq\SX$, let the terminal cost be given by $V_\mathrm{f}\equiv V^\mathrm{ces}_\infty$ and let Assumption~\ref{ass:cont_comp} hold. If the limit $V^\mathrm{ces}_\infty(x)=\lim_{N\to\infty}V^\mathrm{ces}_N(x)$ exists for all $x\in\SX_\mathrm{f}$, then the terminal cost $V_\mathrm{f}$ satisfies Assumption~\ref{ass:terminal_cost_on_orbit}.
	}
	\begin{proof}
		First note that $\ell_{\varPi^\star}=0$ is implied by convergence of the Cesàro cost and that $V^\mathrm{ces}_\infty(x)\leq\ell(x,u)+V^\mathrm{ces}_\infty(f(x,u))$ for all $(x,u)\in\SX_\mathrm{f}\times\SU$ by \eqref{eq:cesaro_iteration}. Hence, for all $i\in\SI_{[0,p^\star-1]}$ we have 
		\begin{align} \label{eq:cesaro_cost_ineq}
			V^\mathrm{ces}_\infty(\varPi_\SX^\star(i))\leq\ell(\varPi^\star(i))+V^\mathrm{ces}_\infty(\varPi^\star_\SX([i+1]_{p^\star})).
		\end{align}
		Now suppose there exists $j\in\SI_{[0,p^\star-1]}$ such that this inequality is strict, i.e., $V^\mathrm{ces}_\infty(\varPi_\SX^\star(j))<\ell(\varPi^\star(j))+V^\mathrm{ces}_\infty(\varPi^\star_\SX([j+1]_{p^\star}))$. Then, applying \eqref{eq:cesaro_cost_ineq} in total $p^\star$ times and exploiting that $[j+p^\star]_{p^\star}=j$ leads to 
		\begin{align*}
			V^\mathrm{ces}_\infty(\varPi_\SX^\star(j))&<V^\mathrm{ces}_\infty(\varPi_\SX^\star([j+p^\star]_{p^\star}))+\!\sum_{i=j}^{j+p^\star-1}\!\!\ell(\varPi^\star([i]_{p^\star})) \\
			& = V^\mathrm{ces}_\infty(\varPi_\SX^\star(j)) +\!\sum_{i=j}^{j+p^\star-1}\!\!\ell(\varPi^\star([i]_{p^\star})).
		\end{align*}
		Since the average cost of $\varPi^\star$ is zero, the sum in the expression above is also zero and hence we arrive at the contradiction $V^\mathrm{ces}_\infty(\varPi_\SX^\star(j))<V^\mathrm{ces}_\infty(\varPi_\SX^\star(j)$, which proves that relation \eqref{eq:cesaro_cost_ineq} must hold with equality for all $i\in\SI_{[0,p^\star-1]}$.
	\end{proof}
	Therefore, under relatively minor assumptions, $V_\mathrm{f}(x)=V^\mathrm{ces}_\infty(x)$ is a suitable terminal cost. Sufficient conditions for the existence of the limit $N\rightarrow\infty$ of $V^\mathrm{ces}_N$ are provided in \cite{Mair25} and include next to strict dissipativity and minimality of $\varPi^\star$ also a local controllability and a finite time reachability condition. Since the infinite-horizon Cesàro cost is in practice hard to compute, we add two brief practical considerations. First, \eqref{eq:cesaro_iteration} provides a simple iterative rule to approximate $V^\mathrm{ces}_\infty(x)$. However, since this is only an approximation, it can only be guaranteed that $\norm{V^\mathrm{ces}_N(x)-V^\mathrm{ces}_\infty(x)}<\eps$ for all $x\in\SX_\mathrm{f}$ and for some $\eps>0$. Therefore, instead of \eqref{eq:terminal_conditions}, only 
	\[
	V_\mathrm{f}(x)\geq V_\mathrm{f}(f(x,u_\mathrm{f}(x))) + \ell(x,u_\mathrm{f}(x))+2\eps
	\]
	is guaranteed, which still guarantees practical asymptotic stability for minimal orbits \cite{Schwenkel2024}, but not the non-averaged performance bound.
	Second, if the optimal periodic orbit $\varPi^\star$ is known, one can set $\SX_\mathrm{f}=\varPi^\star(\SI_{[0,p^\star-1]})$ and determine a terminal cost function such that $V_\mathrm{f}(x)+c=V^\mathrm{ces}_\infty(x)$ for all $x\in\varPi^\star_\SX(\SI_{[0,p^\star-1]})$ with $c\in\R$, i.e. such that $V_\mathrm{f}(\cdot)$ and $V^\mathrm{ces}_\infty(\cdot)$ only differ by a constant. This is done by setting $V_\mathrm{f}(\varPi^\star(0))=0$ and exploiting the relation $V^\mathrm{ces}_\infty(\varPi^\star_\SX(i+1))=V^\mathrm{ces}_\infty(\varPi_\SX^\star(i))-\ell(\varPi^\star(i))$ for $i\in\SI_{[0,p^\star-1]}$. Since the constant $c$ does not change the minimizer of \eqref{eq:finite_horizon_value_function}, this is a valid terminal cost.
	

	\section{Conclusion} \label{sec:conclusion}
	
	We considered an EMPC scheme with terminal conditions for systems with optimal periodic operation, for which we have shown asymptotic stability guarantees for general non-minimal orbits and a non-averaged performance bound for minimal orbits. This extends the works that have derived similar guarantees for optimal operation at steady state. At last, we provided the insight that a suitable terminal cost can be designed by computing an infinite-horizon Cesàro cost.

	\bibliographystyle{IEEEtran}
	\bibliography{../../../literature/Literature}
	
	\appendix
	
	\lemma[Monotonic \!increase \!of \!$\t J^\mathrm{ces}_N$\!]{\label{lem:monotonic_increase_Juc}\!
		\hspace{-0.2cm}\mbox{Let~\!Assumption~\!\ref{ass:strict_dissipativity}} hold and let $x\in\SX_\infty$ and $u\in\SU^\infty(x)$. Then, $\t J^\mathrm{ces}_{N+1}(x,u)\geq\t J^\mathrm{ces}_{N}(x,u)$ and $\t J^\mathrm{ces}_\infty(x,u)\geq\t J^\mathrm{ces}_{N}(x,u)$ for all $N\in\N$. 
	}
	\begin{proof}
		Since $(1-\tfrac{k}{N+1})\geq(1-\tfrac{k}{N})$ for all $k\in\SI_{[0,N]}$ and $\t\ell(x,u)\geq0$, for any $x$ and $u$, we infer $\t J^\mathrm{ces}_{N+1}(x,u)\geq\t J^\mathrm{ces}_{N}(x,u)$. The second inequality follows by induction.
	\end{proof}
	
	\lemma[Bound propagation]{\label{lem:bound_propagation}
		Let Assumptions~\ref{ass:cont_comp}, \ref{ass:strict_dissipativity}, \ref{ass:terminal_conditions} and \ref{ass:bound_V_N} hold and let $\ell_{\varPi^\star}=0$. Then, for any $x\in\SX_\infty$, any $u\in\SU^\infty(x)$ and for each $q\in\N$, the inequality $J^\mathrm{ces}_\infty(x,u)\leq V^\mathrm{uc}_\infty(x)+1$ implies the inequality $J^\mathrm{ces}_\infty(x_u(q,x),u(q+\cdot))\leq V^\mathrm{uc}_\infty(x_u(q,x))+1$.
	}
	\begin{proof}
		Let $x\in\SX_\infty$ and $u\in\SU^\infty(x)$ such that $J^\mathrm{ces}_\infty(x,u)\leq V^\mathrm{uc}_\infty(x)+1$.
		By Lemma~\ref{lem:bounds_on_Vuc} and compactness of $\SX$ it follows that $J^\mathrm{ces}_\infty(x_u(1,x),u(1+\cdot))$ is bounded. Hence, Lemma~\ref{lem:separation_infinite-horizon_costs} implies for any $q\in\N$ that ${J^\mathrm{ces}_\infty(x,u)=J^1_q(x,u)+J^\mathrm{ces}_\infty(x_u(q,x),u(q+\cdot))}$. 
		Then, $J^\mathrm{ces}_\infty(x,u)\leq V^\mathrm{uc}_\infty(x)+1\leq J^1_q(x,u)+V^\mathrm{uc}_\infty(x_u(q,x))+1$ shows the desired result.
	\end{proof}

	\lemma[Separation of infinite-horizon costs]{\label{lem:separation_infinite-horizon_costs}
		Let Assumptions~\ref{ass:cont_comp} and \ref{ass:strict_dissipativity} hold, let $q\in\N$, $x\in\SX_\infty$ and $\ell_{\varPi^\star}=0$. Then, for any $C>0$ and $u\in\SU^\infty(x)$ satisfying $J^\mathrm{ces}_\infty(x,u)\leq C$, it holds that ${J^\mathrm{ces}_\infty(x,u)=J^1_q(x,u)+J^\mathrm{ces}_\infty(x_u(q,x),u(q+\cdot))}$.
	}	
	\begin{proof}
		A calculation \cite{Mair25} yields for any $N\in\N$
		$
		J_N^\mathrm{uc}(x,u)=\ell(x,u(0))+(1-\frac{1}{N})J^\mathrm{ces}_{N-1}(x_u(1,x),u(1+\cdot)).
		$
		Taking the $\limsup_{N\rightarrow\infty}$ of the expression above results in $J^\mathrm{ces}_\infty(x,u)=\ell(x,u(0))+J^\mathrm{ces}_\infty(x_u(1,x),u(1+\cdot))$, since $J^\mathrm{ces}_\infty(x_u(1,x),u(1+\cdot))$ is upper bounded by assumption and lower bounded by Lemma~\ref{lem:bounds_on_Vuc}. Applying the steps above $q$ times proves the claim.
	\end{proof}


\end{document}